\documentclass[12pt,twoside]{article}

\usepackage{graphicx}

\usepackage{amsthm}

\usepackage{amsmath}

\usepackage{amsfonts}

\usepackage{amssymb}

\usepackage[all]{xy}

\newtheorem{theorem}{Theorem}

\newtheorem{definition}{Definition}

\newtheorem{example}{Example}

\newtheorem{proposition}{Proposition}

\newtheorem{remark}{Remark}

\date{} 

\begin{document} 

\centerline{\bf Int. J. Contemp. Math. Sciences, Vol. 2, 2007, no. 26, 1293 - 1305} 

\centerline{} 

\centerline{} 

\centerline {\Large{\bf On the variety of two dimensional}} 

\centerline{} 

\centerline{\Large{\bf real associative algebras}} 

\centerline{} 

\centerline{\bf {Jos\'{e} Mar\'{\i}a Ancochea Berm\'{u}dez}\footnote[1]{This author has
    been partially supported by the research proyect MTM2006-09152 of
    Ministerio de Educaci\'{o}n y Ciencia}}
\centerline{} 

\centerline{Dpto. Geometr\'{\i}a y Topolog\'{\i}a. Facultad de
 Matem\'{a}ticas} 

\centerline{Universidad Complutense de Madrid} 

\centerline{Plaza de Ciencias, 3 28040 Madrid, Spain} 

\centerline{ancochea@mat.ucm.es} 

\centerline{} 

\centerline{\bf {Javier Fres\'{a}n}} 

\centerline{} 

\centerline{Dpto. Geometr\'{\i}a y Topolog\'{\i}a. Facultad de
 Matem\'{a}ticas} 

\centerline{Universidad Complutense de Madrid} 

\centerline{Plaza de Ciencias, 3 28040 Madrid, Spain} 

\centerline{jfresan@estumail.ucm.es} 

\centerline{} 

\centerline{\bf {Jonathan S\'{a}nchez Hern\'{a}ndez}} 

\centerline{} 

\centerline{Dpto. Geometr\'{\i}a y Topolog\'{\i}a. Facultad de
 Matem\'{a}ticas} 

\centerline{Universidad Complutense de Madrid} 

\centerline{Plaza de Ciencias, 3 28040 Madrid, Spain} 

\centerline{jnsanchez@mat.ucm.es} 

\newtheorem{Theorem}{\quad Theorem}[section] 

\newtheorem{Definition}[Theorem]{\quad Definition} 

\newtheorem{Corollary}[Theorem]{\quad Corollary} 

\newtheorem{Lemma}[Theorem]{\quad Lemma} 

\newtheorem{Example}[Theorem]{\quad Example} 

\begin{abstract} This paper consists of a description of the variety of two 
dimensional associative algebras within the framework of Nonstandard Analysis. 
By decomposing each algebra in $A^{2}$ as sum of a Jordan algebra and a 
Lie algebra, we calculate the isomorphism classes of two dimensional 
associative algebras over the field of real numbers and determine the 
open components and the contractions of the variety.
\end{abstract} 



\section{Introduction} 

The aim of this work is to study some properties of the variety of two 
dimensional associative algebras, specially those concerning rigidity and 
contractions. We first obtain a decomposition of associative algebras as 
sum of a Jordan algebra and a Lie algebra, which enables us to use known 
results on Jordan algebras \cite{AnCa} to classify two dimensional 
associative algebras over the field of real numbers. Then we introduce the 
concept of perturbation within Nelson's Internal Set Theory \cite{Ne} and 
derive the perturbation equations of the variety. Nonstandard Analysis tools 
permit us to prove that $\mathcal{A}^2$ has four open components, two of dimension $4$ 
and two of dimension $2$. The remaining algebras of the variety are obtained by 
contraction of the rigid algebras which define the open components.

\begin{definition}
An associative algebra law over $\mathbb{R}$ is a bilinear mapping \newline 
$\beta:\mathbb{R}^{n} \times \mathbb{R}^{n}\longrightarrow \mathbb{R}^{n}$ 
satisfying the constraint \begin{equation}
\beta(\beta(x,y),z)-\beta(x,\beta(y,z))=0, \quad x,y, z \in \mathbb{R}^{n}. 
\end{equation} We will abbreviate $\beta(x, y)$ by $x\circ 
y,$ and $\mathcal{A}^{n}$ will denote the set of associative algebras over 
$\mathbb{R}^{n}.$ The ordered pair $(\mathbb{R}^{n},\beta)$ is called associative algebra. 
\end{definition}

\begin{definition}
A Jordan algebra law over $\mathbb{R}$ is a symmetric bilinear mapping 
$\varphi: \mathbb{R}^{n} \times \mathbb{R}^{n} \longrightarrow \mathbb{R}^{n}$ 
which verifies the identity
\begin{equation}
\varphi(\varphi(x,x),
\varphi(x,y))-\varphi(x,\varphi(\varphi(x,x),y))=0, \quad x, y \in \mathbb{R}^{n}.  \nonumber
\end{equation} The ordered pair $(\mathbb{R}^{n}, \varphi)$ is a Jordan 
algebra and $\mathcal{J}^n$ will denote the set of Jordan algebra laws.
\end{definition}

\begin{definition}
A Lie algebra law over $\mathbb{R}^{n}$ is an alternate bilinear mapping 
$\mu: \mathbb{R}^{n} \times \mathbb{R}^{n} \longrightarrow \mathbb{R}^{n}$ 
satisfying the Jacobi identity
\begin{equation}
\mu(\mu(x,y), z)+\mu(\mu(y,z),
x)+\mu(\mu(z,x), y)=0, \quad x, y, z \in \mathbb{R}^{n}. \nonumber
\end{equation} We will denote the set of Lie algebras over $\mathbb{R}^{n}$ by $\mathcal{L}^n.$ The ordered pair $(\mathbb{R}^{n}, \mu),$ where $\mu \in \mathcal{L}^{n},$ is a Lie algebra. 
\end{definition}

Now we enunciate a proposition which plays a fundamental role in the 
classification which will be given.

\begin{proposition}
Let $(\mathbb{R}^{n}, \circ)$ be a real associative algebra. Then:
\begin{enumerate}
\item The law $\varphi$ defined by \begin{equation} \varphi(x,
y)=\frac{x \circ y + y \circ x} {2},  \quad x, y \in \mathbb{R}^{n} \nonumber
\end{equation} is a Jordan algebra law.
\item The law $\mu$ defined by \begin{equation} \mu(x,
y)=\frac{x \circ y - y \circ x} {2}, \quad x, y \in \mathbb{R}^{n} \nonumber
\end{equation} is a Lie algebra law. \end{enumerate} \end{proposition}

Let us recall that if there exists a non-zero vector $u$ such that 
$\beta(u, v)=0$ for all $v \in \mathbb{R}^{n},$ we say 
that $\beta$ has left isotropy. Similarly, if $\beta(v, 
u)=0$ for all $v \in \mathbb{R}^{n},$ then $u$ is a right 
isotropic vector. An associative algebra is said to be simple if it does not 
admit any proper ideal.

Let $B=\{e_1,\ldots,e_n\}$ be a basis for $\mathbb{R}^{n}.$ 
It is possible to identify each algebra in $\mathcal{A}^n$ with its structure 
constants, that is, to consider $\beta \in \mathcal{A}^n$ as the tensor $(a_{ij}^{k}) \in \mathbb{R}^{n^{3}}$ whose coordinates, univocally determined by $e_{i} \circ e_{j}=a_{ij}^{k}e_{k},$ are the solutions of the system
\begin{equation}
a_{ij}^{h}a_{hk}^{l}=a_{ih}^{l}a_{jk}^{h}, \quad 1 \leq i,j,k,h,l
\leq n. \label{sys} \end{equation} 
That gives $\mathcal{A}^n$ a structure of algebraic variety embedded in 
$\mathbb{R}^{n^{3}}.$ From now on we will identify each algebra with its law. In the case of associative algebras, non-written products will be supossed to be zero. If $\varphi$ is a Jordan algebra or a Lie algebra only non-zero products $\varphi(e_i,e_j),$ with $i \leq j,$ will be written.   

\section{Classification of two dimensional real associative algebras}

If $n=2,$ real associative algebras are given by the relations
\begin{align*}
e_{1} \circ e_{1} &= a_{1} e_{1} + a_{2}e_{2},
\\ e_{1} \circ e_{2} &= b_{1} e_{1} +
b_{2}e_{2}, \\ e_{2} \circ e_{1} &= c_{1} e_{1}
+ c_{2}e_{2}, \\ e_{2} \circ e_{2} &= d_{1}
e_{1} + d_{2}e_{2}.
\end{align*} Or equivalently by a coefficient matrix of the form:
\[
\left(\begin{array}{cc} a_{1} & a_{2} \\ b_{1} & b_{2} \\ c_{1} &
c_{2} \\ d_{1} & d_{2} \\
\end{array} \right)
\]

Developing \eqref{sys} the following equations are obtained:
\begin{align}
\begin{split}
a_{2}b_{1}&=a_{2}c_{1}, \\ 
a_{2}b_{2}&=a_{2}c_{2}, \\ 
b_{1}b_{2}&=a_{2}d_{1}, \\ 
a_{2}d_{1}&=c_{1}c_{2}, \\ 
a_{2}b_{1}+b_{2}^{2}&=a_{1}b_{2}+a_{2}d_{2}, \\ 
a_{1}c_{1}+b_{1}c_{2}&=a_{1}b_{1}+b_{2}c_{1}, \\ 
a_{1}d_{1}+b_{1}d_{2}&=b_{1}^{2}+b_{2}d_{1}, \\ 
a_{1}c_{2}+a_{2}d_{2}&=a_{2}c_{1}+c_{2}^{2}, \\ 
b_{1}c_{2}+b_{2}d_{2}&=b_{2}c_{1}+c_{2}d_{2}, \\ 
c_{1}^{2}+c_{2}d_{1}&=a_{1}d_{1}+c_{1}d_{2}. \label{ecua}
\end{split}
\end{align} Thus, $A^{2}$ is an algebraic variety embedded in 
$\mathbb{R}^{8}$ and defined by the above system of homogeneous polynomials.

\smallskip

Let us consider the natural action of the general linear group $GL(n,\mathbb{R})$ over the variety: \[\begin{array}{rcl}
GL(n,\mathbb{R}) \times \mathcal{A}^{n}& \longrightarrow &\mathcal{A}^{n} \\ (f, \beta)
&\longmapsto & f^{-1} \ast \beta \ast (f \times f),
\end{array}\] with  $f^{-1} \ast \beta \ast (f \times f) (x, 
y)= f^{-1}(\beta(f(x),f(y)).$ For each $\beta \in
\mathcal{A}^{n},$ the orbit under this action, $\mathcal{O}(\beta),$ represents the set of associative algebras isomorphic to $\beta.$ Our first aim is to determine the space of orbits
\[
A^{2}/GL(2,\mathbb{R})=\{\mathcal{O}(\beta)\}_{\beta \in \mathcal{A}^{2}},
\]
that is, the isomorphism classes of two dimensional real 
associative algebras. 

For that purpose, we will make use of the decomposition
\[
x \circ y = \frac{x \circ y+y \circ
x}{2}+\frac{x \circ y-y \circ
x}{2}=\varphi(x, y) + \mu(x, y),
\] where $\varphi$ is a Jordan algebra and $\mu$ is a Lie 
algebra defined by the product $\mu(e_1,e_2)=ae_1+be_2.$ In \cite{AnCa} a classification theorem for $\mathcal{J}^2$ is proved:

\begin{theorem} \label{teor1}
Let $\varphi$ be a non-Abelian two dimensional real Jordan algebra. Then $\varphi$ is isomorphic to one of the following pairwise non-isomorphic Jordan 
algebras:
\begin{enumerate}
\item $\varphi_{1}(e_1,e_1)=e_1,
\ \varphi_{1}(e_1,e_2)=e_2, \ %
\varphi_{1}(e_2,e_2)=-e_1.$

\item $\varphi_{2}(e_1,e_1)=e_1,
  \  \varphi_{2}(e_1,e_2)=e_2, \ 
  \varphi_{2}(e_2,e_2)=e_1.$

\item $\varphi_{3}(e_1,e_1)=e_1,
\  \varphi_{3}(e_1,e_2)=e_2, \ \
\varphi_{3}(e_2,e_2)=0.$

\item $\varphi_{4}(e_1,e_1)=0, \ \
\varphi_{4}(e_1,e_2)=0, \ \ \
\varphi_{4}(e_2,e_2)=e_2.$

\item $\varphi_{5}(e_1,e_1)=e_2,
\  \varphi_{5}(e_1,e_2)=0, \ \ \ 
\varphi_{5}(e_2,e_2)=0.$

\item $\varphi_{6}(e_1,e_1)=e_1,
\  \varphi_{6}(e_1,e_2)=\frac{1}{2}e_2, \ 
\varphi_{6}(e_2,e_2)=0.$
\end{enumerate}
\end{theorem}

\begin{example} \label{ex1}
To illustrate how this result may be applied to solve the classification 
problem, let us consider the case in which $\varphi$ is isomorphic to $\varphi_1.$ Then, there exists a basis $\{e_1,e_2\}$ such that:
\begin{align*}
e_1 \circ
e_1&=\varphi_{1}(e_1,e_1)+\mu(e_1,e_1) 
=e_1. \\ e_1 \circ e_2& =
\varphi_{1}(e_1, e_2)+\mu(e_1,e_2) 
=ae_1+(1+b)e_2. \\
e_2 \circ e_1&=\varphi_{1}(e_2,
e_2)+\mu(e_2,e_1)=-ae_1+(1-b)e_2. \\ 
e_2 \circ
e_2&=\varphi_{1}(e_2,
e_2)+\mu(e_2,e_2)=-e_1.
\end{align*}
Structure constants must satisfy \eqref{ecua}, so $a=b=0$ and $\beta$ coincides with $\varphi_1.$
\end{example}

Making the same calculus for each $\varphi_i \in \mathcal{J}^2,$ the following theorem is obtained. For the first five Jordan algebras we have $a=b=0,$ but when we consider $\varphi$ isomorphic to $\varphi_6,$ the system \eqref{ecua} admits two different solutions: $(0,\frac{1}{2})$ and $(0,-\frac{1}{2}).$ Thus, there are seven isomorphism classes in $\mathcal{A}^2.$

\begin{theorem} \label{teor2}
Let $\beta$ be a two dimensional real associative algebra. If $\beta$ is not Abelian, then $\beta$ is isomorphic to one of the following pairwise non-isomorphic associative algebras:
\begin{enumerate}
\item $\beta_{1}: \  e_1 \circ e_1=e_1, \ 
e_1 \circ e_2=e_2, \  e_2 \circ
{e}_1=e_2, \  e_2 \circ e_2=-e_1.$

\item $\beta_{2}: \  e_1 \circ e_1=e_1, \ 
e_1 \circ e_2=e_2, \  e_2 \circ
e_1=e_2, \  e_2 \circ e_2=e_1.$

\item $\beta_{3}: \  e_1 \circ e_1=e_1, \ 
e_1 \circ e_2=e_2, \  e_2 \circ
e_1=e_2, \  e_2 \circ e_2=0.$

\item $\beta_{4}: \  e_1 \circ e_1=0, \ \ 
e_1 \circ e_2=0, \ \  e_2 \circ
e_1=0, \ \  e_2 \circ e_2=e_2.$

\item $\beta_{5}: \  e_1 \circ e_1=e_2, \ 
e_1 \circ e_2=0, \ \ e_2 \circ
e_1=0, \ \  e_2 \circ e_2=0.$

\item $\beta_{6}: \  e_1 \circ e_1=e_1, \ 
e_1 \circ e_2=e_2, \  e_2 \circ
e_1=0, \ \  e_2 \circ e_2=0.$

\item $\beta_{7}: \  e_1 \circ e_1=e_1, \ 
e_1 \circ e_2=0,\ \  e_2 \circ
e_1=e_2,\  e_2 \circ e_2=0.$
\end{enumerate}
Moreover, $\beta_2$ is the unique simple algebra in $\mathcal{A}^2.$
\end{theorem}

\begin{remark}
Considered over $\mathbb{C},$ $\beta_{1}$ and $\beta_{2}$ are isomorphic, 
with the change of basis given by $x_{1}=e_1$ and 
$x_{2}= i e_2.$ Thus, we have found an example of a simple 
associative algebra whose complexification is not simple. In $\mathcal{A}^2,$ 
$\beta_{3}, \beta_{4}$ and $\beta_{5}$ have right and left isotropy, 
$\beta_{6}$ has just left isotropy and $\beta_7$ has just right isotropy.
\end{remark}

\begin{proposition} \label{prop2}
Let $\beta \in \mathcal{A}^2$ be a two dimensional real associative algebra.
\begin{enumerate}
\item If $\beta$ has no isotropy, then $\beta$ is isomorphic either to 
$\beta_1$ or to $\beta_2.$
\item If $\beta$ has just left isotropy, then $\beta$ is isomorphic to 
$\beta_{6}.$
\item If $\beta$ has just right isotropy, then $\beta$ is isomorphic to 
$\beta_{7}.$
\end{enumerate}
\end{proposition}

\section{Perturbations of associative algebras}
Considering in $\mathcal{A}^n$ the subspace topology induced by $\mathbb{R}^{n^{3}},$ 
it is possible to give the following definition:

\begin{definition}
An associative algebra $\beta$ is rigid in $\mathcal{A}^n$ if its orbit under the action of $GL(n,\mathbb{R})$ is open.  
\end{definition}

To study rigidity in the framework of Internal Set Theory, the concept of 
perturbation is introduced (cf. \cite{Di}, \cite{LuGo}, \cite{Ne} for the details).

\begin{definition}
Let $n$ be standard and let $\beta_{0}$ be a real associative algebra. We say that $\beta$ is a perturbation of $\beta_{0}$ if $\beta(x,y)$ and $\beta_{0}(x,y)$ are infinitely close for all $x, y$ standard vectors of $\mathbb{R}^{n}.$
\end{definition}

\begin{proposition}
A standard associative algebra law $\beta_{0} \in \mathcal{A}^{n}$ is rigid if and only if any perturbation of $\beta_{0}$ is isomorphic to it.
\end{proposition}

\begin{proof}
If $\beta_{0} \in \mathcal{A}^{n}$ is rigid, then $\mathcal{O}(\beta_{0})$ is open. Thus, $\mathcal{O}(\beta_{0})$ contains the halo of $\beta_{0}$ and any perturbation $\beta$ of $\beta_{0}$ is isomorphic to $\beta_{0}.$ Conversely, if any perturbation belongs to $\mathcal{O}(\beta_{0}),$ then the halo of $\beta_{0}$ is contained in $\mathcal{O}(\beta_{0})$ and, by transference, 
$\mathcal{O}(\beta_{0})$ is open, so $\beta_{0}$ is rigid.
\end{proof}

Now we enunciate a theorem, due to M. Goze~\cite{Go}, basic to determine 
the rigidity of associative algebras.

\begin{theorem}
Let $M_0$ be a standard point in $\mathbb{R}^{n}.$ Then every point $M \in 
\mathbb{R}^{n}$ infinitely close to $M_0$ admits a decomposition of the 
form:
\begin{equation}
M=M_{0}+\epsilon_{1} v_{1}+\epsilon_{1} \epsilon_{2}
v_{2}+\ldots +\epsilon_{1} \epsilon_{2} \ldots \epsilon_{p}
v_{p},\nonumber
\end{equation} where $\epsilon_{1},\epsilon_{2},\ldots, \epsilon_{p}$ are infinitely small scalars and $v_{1},v_{2}, \ldots, v_{p}$ 
linearly independent vectors. Moreover, if $M=M_{0}+\eta_{1} u_{1}+\eta_{1} \eta_{2}u_{2}+\ldots +\eta_{1} \eta_{2} \ldots \eta_{q} u_{q}$ is 
another decomposition of $M,$ then $p=q$ and the flag defined by 
$v_{1},v_{2}, \ldots, v_{p}$ coincides with the flag defined 
by $u_{1},u_{2}, \ldots, u_{q}.$ The integer $p$ which 
describes the equivalence class of a point is called length of $M.$
\end{theorem}

As a consequence of Goze's theorem, any perturbation of the standard 
law $\beta_{0} \in \mathcal{A}^{n}$ may be written as
\begin{equation}
\beta=\beta_{0}+\epsilon_{1} \varphi_{1}+\epsilon_{1} \epsilon_{2}
\varphi_{2}+\ldots +\epsilon_{1} \epsilon_{2} \ldots \epsilon_{p}
\varphi_{p}, \label{goz}
\end{equation} where $\epsilon_{i}$ are infinitely small and $\varphi_{i}: 
\mathbb{R}^{n} \times \mathbb{R}^{n} \rightarrow \mathbb{R}^{n}$ are 
independent bilinear mappings. Considering $\beta_0$ a standard element of 
the variety and $\beta$ a perturbation of $\beta_0$ decomposed according to 
\eqref{goz}, we have that the shade of the straight line which joins $\beta$ and 
$\beta_0$ is a standard straight line which belongs to the tangent cone to 
the variety in $\beta_0.$

\medskip

Let $\beta_{0} \in \mathcal{A}^{n}$ be a standard associative algebra law and let $f \in GL(n,\mathbb{R})$ be standard. For any $\epsilon$ infinitely small, 
the endomorphism $Id+\epsilon f$ belongs to the general linear group, so it 
makes sense to consider the action
\[
(Id+\epsilon f)^{-1}\beta_{0}((Id+\epsilon f),(Id+\epsilon
f))=\beta_{0}+\epsilon(\delta_{\beta_{0}}f)+\epsilon^{2}(\Delta(\beta_{0},f,\epsilon)),
\]
where 
\[\delta_{\beta}f(x,y)=\beta(f(x),
y)+\beta(x,f(y))-f(\beta(x,y))
\] are the $2-$coboundaries of Hochschild cohomology. Moreover, since the straight line which joins $\beta_0$ and an 
infinitely close point $\beta_{0}'$ is tangent to the orbit 
$\mathcal{O}(\beta_{0})$ in $\beta_{0},$ the tangent space is given by
\begin{equation} T_{\beta_{0}}\mathcal{O}(\beta_{0})=\{\delta_{\beta_{0}}f 
\,:\, f \in GL(n, \mathbb{R})\}. \label{tg} \end{equation} The dimension of the orbit is the 
vector dimension of $T_{\beta_{0}}\mathcal{O}(\beta_{0}).$

\medskip

Now let $\beta_{0} \in \mathcal{A}^{n}$ be a standard associative algebra and let us consider a perturbation $\beta$ of $\beta_{0}.$ According to Goze's theorem, for an integer $p \leq n,$ there exist $\epsilon_{1}, \ldots, \epsilon_{p}$ infinitely small, with $\epsilon_{1} \neq 0$, and $\varphi_{1}, \ldots, \varphi_{p}$ independent bilinear mappings such that:
\[
\beta=\beta_{0}+\epsilon_{1} \varphi_{1}+\epsilon_{1} \epsilon_{2}
\varphi_{2}+ \ldots + \epsilon_{1}\ldots \epsilon_{p} \varphi_{p}.
\] Denoting by $\beta_{1} \circ \beta_{2}$ the trilinear mapping 
defined as \[
\beta_{1} \circ \beta_{2}(x, y,
z)=\beta_{1}(\beta_{2}(x, y),
z)-\beta_{1}(x, \beta_{2}(y,z))
+\beta_{2}(\beta_{1}(x,y),z) -
\beta_{2}(x,\beta_{1}(y,z)), \nonumber
\] $\beta$ is an associative algebra if and only if $\beta \circ 
\beta \equiv 0.$ If the infinitesimal part is denoted by $\xi,$ then 
$\beta=\beta_{0}+\xi$ and the latter condition is written:
\begin{equation*}
\beta \in A^{n} \Leftrightarrow (\beta_{0}+\xi) \circ (\beta_{0}+\xi)
= \beta_{0} \circ \beta_{0}+\beta_{0} \circ \xi + \xi \circ \beta_{0}
+ \xi \circ \xi =0,
\end{equation*} Since $\beta_{0} \in \mathcal{A}^{n}$ and $\beta_{0} \circ \xi = \xi 
\circ \beta_{0}=\delta^2_{\beta_0}\xi,$ we have: \begin{equation}
2\delta_{\beta_{0}}^{2}\xi + \xi \circ \xi =0,
\end{equation} where \begin{equation*}
\delta_{\beta}^{2}\varphi(x, y, 
z)=\beta(\varphi(x, y), 
z)-\beta(x,\beta(y,z))+\varphi(\beta(x,y),z) 
-\varphi(x,\beta(y,z)) \end{equation*} are the $2-$cocycles of Hochschild 
cohomology. 

This is the perturbation equation of $\beta_{0}.$ Developing the 
expression, dividing by $\epsilon_{1},$ and separating the standard and the 
infinitesimal part, we obtain:\begin{align*} 
\delta_{\beta_{0}}\varphi_{1}&=0, 
\\
2(\epsilon_{2}\delta_{\beta_{0}}^2\varphi_{2}+ \ldots + \epsilon_{2}\ldots
\epsilon_{p}\delta_{\beta_{0}}^2\varphi_{p} +\epsilon_{1} \varphi_{1} \circ 
\varphi_{1} +
\epsilon_{1} \epsilon_{2}^{2}\varphi_{2} \circ \varphi_{2}+ & \\
+\ldots +\epsilon_{1} \epsilon_{2}^{2} \ldots \epsilon_{p}^{2}\varphi_{p} \circ
\varphi_{p})+\epsilon_{1} \epsilon_{2} \varphi_{1} \circ
\varphi_{2}+\ldots+\epsilon_{1} \ldots \epsilon_{p} \varphi_{1} \circ
\varphi_{p}+ &\\
+\epsilon_{1}\epsilon_{2}^{2}\epsilon_{3}\varphi_{2} \circ
\varphi_{3}+ \ldots + \epsilon_{1}\epsilon_{2}^{2} \ldots
\epsilon_{p-1}^{2}\epsilon_{p}\varphi_{p-1} \circ \varphi_{P}&=0.\end{align*}

\section{Rigid laws}
The interest of determining which associative algebras are rigid lies in the 
main role they play in the study of the variety, because their orbits are 
the open components of $\mathcal{A}^n.$ In this paragraph we calculate
the dimension of the orbits of the associative algebra laws obtained in the 
classification theorem and prove which of them are rigid. From now 
on, $\varphi_i$ will denote the Jordan algebras of theorem \ref{teor1}, and $\beta_i$ the associative algebras of theorem \ref{teor2}.

\medskip

Let $f \in GL(2,\mathbb{R}),$ with 
$f(e_{1})=ae_{1}+be_{2}$ and 
$f(e_{2})=ce_{1}+de_{2}$, be a non-singular endomorphism. 
Evaluating \eqref{tg} for each $\beta_{i},$ the following tangent spaces to the 
orbits are obtained:

\medskip
\begin{center}
\begin{tabular}{|c|}\hline
$T_{\beta_{1}}f=\left(\begin{array}{cc} a & b \\ -b & a \\ -b & a \\
a-2d & b+2c \\
\end{array} \right) \quad T_{\beta_{2}}f=\left(\begin{array}{cc}
a & b \\ b & a \\ b & a \\ 2d-a & 2c-b \\
\end{array} \right)$
\\\hline $ T_{\beta_{3}}f=\left(\begin{array}{cc} a & b \\ 0 & a \\ 0
& a \\ 0 & c \\
\end{array} \right) \quad T_{\beta_{4}}f=\left(\begin{array}{cc}
0 & 0 \\ 0 & b+d \\ 0 & b+d \\ -c & d \\
\end{array} \right)$
\\\hline $T_{\beta_{5}}f=\left(\begin{array}{cc} -c & 2a-d \\ 0 & c \\
0 & c \\ 0 & 0 \\
\end{array} \right) \quad T_{\beta_{6}}f=\left(\begin{array}{cc}
a & 0 \\ 0 & a \\ c & 0 \\ 0 & c \\
\end{array} \right) \quad T_{\beta_{7}}f=\left(\begin{array}{cc}
a & 0 \\ c & 0 \\ 0 & a \\ 0 & c \\
\end{array} \right)$ \\\cline{1-1}%
\end{tabular}
\end{center}
\medskip

Therefore, the dimension of the orbits are:
\begin{gather*}
\dim \mathcal{O}(\beta_{1})=\dim \mathcal{O}(\beta_{2})=4, \\
\dim \mathcal{O}(\beta_{3})=\dim \mathcal{O}(\beta_{4})=3, \\
\dim \mathcal{O}(\beta_{5})=\dim \mathcal{O}(\beta_{6})=\dim \mathcal{O}(\beta_{7})=2.
\end{gather*}

To determine the open components of $\mathcal{A}^2,$ we will apply theorem \ref{teor2} and some properties of the variety of Jordan algebras. In particular, we will make use of the following theorem~\cite{AnCa}:

\begin{theorem}
The only rigid algebras in $\mathcal{J}^{2}$ are $\varphi_{1},$ $\varphi_{2}$ and $\varphi_{6}.$
\end{theorem}

\begin{theorem}
The only two dimensional real associative algebras which are rigid are 
$\beta_{1},$ $\beta_{2},$ $\beta_{6}$ and $\beta_{7}.$ Thus, $\mathcal{A}^2$ has two 
open components of dimension 4 and two open components of dimension 2.
\end{theorem}

\begin{proof}
Let $\beta_{1}=\varphi_{1}+\mu_{1}$ be, with $\varphi_{1}$ and $\beta_{1}$ the 
associated Jordan and Lie algebras respectively. If $\beta$ is a 
perturbation of $\beta_{1},$ then $\beta$ admits a decomposition of the form 
$\beta=\varphi+\mu,$ where $\varphi \sim \varphi_{1}$ and $\mu \sim \mu_{1}.$ Since 
$\varphi_{1}$ is rigid, there exists a basis $\{e_{1}, e_{2}\}$ 
such that $\varphi=\varphi_{1}.$ Then, via the classification theorem (see example \ref{ex1}), $\mu \equiv 0$ and $\beta$ is isomorphic to $\beta_{1}.$ An analogous reasoning proves the rigidity of $\beta_{2}.$

To prove that $\beta_6$ is rigid, let us consider a perturbation 
$\beta=\varphi+\mu,$ with $\varphi \sim \varphi_{6}.$ Since $\varphi_6$ is rigid, $\varphi$ is isomorphic to $\varphi_{6}$ and it is possible to find an infinitely close basis $\{x_1, x_2\},$ where the structure constants 
of $\varphi$ coincide with those of $\varphi_{6}$ and there exist $\epsilon$ and $\epsilon '$ infinitely small such that
\begin{equation*}
\mu(x_1,x_2)=\epsilon x_1+(\frac{1}{2}+\epsilon ')x_2, 
\end{equation*} According to our classification, the associativity condition imposes that $\epsilon=0$ and $\frac{1}{2}+\epsilon '=\pm \frac{1}{2}.$ Since 
$\frac{1}{2}$ and $-\frac{1}{2}$ are not infinitely close, $\epsilon '=0$ 
and $\beta$ is isomorphic to $\beta_{6}.$ The same reasoning proves the 
rigidity of the algebra $\beta_{7}.$

Now let $\beta \in \mathcal{A}^{2}$ be a standard law non-isomorphic to any of the latter algebras. Then there exists a basis $\{e_{1},e_{2}\}$ in 
which the structure constants of $\beta$ are given by one of the following 
matrices:
\begin{equation}
\left(\begin{array}{cc} 1 & 0 \\ 0 & 1 \\ 0 & 1 \\ 0 & 0 \\
\end{array} \right), \left(\begin{array}{cc}
0 & 0 \\ 0 & 0 \\ 0 & 0 \\ 0 & 1 \\
\end{array} \right), \left(\begin{array}{cc}
0 & 1 \\ 0 & 0 \\ 0 & 0 \\ 0 & 0 \\
\end{array} \right). \nonumber
\end{equation} We may consider the perturbations \begin{equation}
\left(\begin{array}{cc} 1 & 0 \\ 0 & 1 \\ 0 & 1 \\ \epsilon & 0 \\
\end{array} \right), \left(\begin{array}{cc}
\epsilon & 0 \\ 0 & 0 \\ 0 & 0 \\ 0 & 1 \\
\end{array} \right), \left(\begin{array}{cc}
0 & 1 \\ 0 & 0 \\ 0 & 0 \\ \epsilon & 0 \\
\end{array} \right), \nonumber
\end{equation} all of them without isotropy. Then, by proposition \ref{prop2}, these perturbations are isomorphic either to $\beta_{1}$ or to $\beta_{2},$ and the laws $\beta_{3},$ $\beta_{4}$ and $\beta_{5}$ are not rigid.\end{proof}

\section{Contractions of associative algebras}
Using the action of the general linear group over the variety of real
associative algebras, it is possible to define a formal notion of
limit in $\mathcal{A}^n$, in analogy with the theory of contractions developed
for Lie~\cite{Go} and Jordan algebras~\cite{AnCa}.

\begin{definition}
Let $\beta_0 \in \mathcal{A}^n$ be a real associative algebra and let $\{f_{t}\} \subset GL(n,\mathbb{R})$ be a family of non-singular endomorphisms depending on a continous parametre $t.$ If the limit
\begin{equation}
\beta(x,y):=\lim_{t\rightarrow0}
f_{t}^{-1}\circ\beta_0(f_{t}(x),f_{t}(y)) \label{lim}
\end{equation} exists for all $x,y \in \mathbb{R}^{n},$ $\beta$ is an associative algebra called contraction of $\beta_0.$
\end{definition}

\begin{example}
$\beta_{3}$ is a contraction of $\beta_{1}$ by the linear transformations \begin{equation}
f_{t}(e_{1})=e_{1}, \quad f_{t}(e_{2})=te_{2}. 
\nonumber \end{equation} Let us consider in $\beta_{1}$ the 
transformed basis $\{x_1=f_{t}(e_{1}),x_2=f_{t}(e_{2})\},$ where vector products are given by:
\begin{align*}
x_1\circ x_1&=e_{1}\circ e_{1}=e_{1}=x_1
\\
x_1\circ x_2&=te_{1}\circ e_{2}=te_{2}=x_2
\\
x_2\circ x_1&=te_{2}\circ e_{1}=te_{2}=x_2
\\
x_2\circ x_2&=t^{2}e_{2}\circ e_{2}=-t^{2}e_{1}=-t^{2}x_1 
\end{align*} Thus, the structure constants of $\beta_{1}$ are represented 
by the matrix \begin{equation*}
\left(\begin{array}{cc} 1 & 0 \\ 0 & 1 \\ 0 & 1 \\ -t^{2} & 0 \\
\end{array} \right).
\end{equation*} It is immediate that (\ref{lim}) holds for $e_1, e_2.$ In the limit, an algebra isomorphic to $\beta_{3}$ is obtained. Therefore, $\beta_3$ is a contraction of $\beta_1.$
\end{example}

\medskip

It is easy to prove from the definition that a contraction of $\beta_0$
corresponds to a closure point of $\mathcal{O}(\beta_0).$ In particular,
rigid algebras are not obtained as a contraction of any non-isormorphic algebra in $\mathcal{A}^n$~\cite{Ni}. It is obvious that the change of basis
$x_i=te_{i},$ $i=1,\ldots, n$ induces a contraction of
any associative algebra over the Abelian algebra. Moreover, for every contraction the following inequality holds:
\begin{equation*}
\dim \mathcal{O}(\beta_0)>\dim \mathcal{O}(\beta)
\end{equation*} That gives us a first criterion to study the contractions of 
the variety. If $\beta$ is a contraction of $\beta_0,$ then its associated
Jordan algebra $\varphi$ is also a contraction of $\varphi_0.$ We have
already proved~\cite{AnCa} that $\beta_4$ is not a contraction of $\beta_1$ in $\mathcal{J}^2.$ The following theorem specifies how to obtain the remaining contractions.

\begin{theorem}
Let $\beta_i$ be the associative algebras of the classification theorem. 
Then $\beta_{3},$ $\beta_{4}$ and $\beta_{5}$ are the only algebras in 
$\mathcal{A}^2$ which appear as the contraction of an associative algebra. More 
precisely:
\begin{enumerate}
\item $\beta_{3}$ is a contraction of $\beta_{1}$ and $\beta_{2}.$
\item $\beta_{4}$ is a contraction of $\beta_{2}.$
\item $\beta_{5}$ is a contraction of $\beta_{1},$ $\beta_{2},$
$\beta_{3}$ and $\beta_{4}$
\end{enumerate}
\end{theorem}

\begin{proof}
Since $\beta_{1},$ $\beta_{2},$ $\beta_{6}$ and $\beta_{7}$ are 
rigid algebras, they are not contractions of any other associative algebra.
\begin{enumerate}
\item We have already proved that when we consider the family of non-singular 
endomorphisms
\begin{equation}
f_{t}(e_{1})=e_{1}, \quad f_{t}(e_{2})=te_{2}, 
\nonumber
\end{equation} a contraction of $\beta_1$ over $\beta_3$ is obtained. The 
same family of linear transformations defines the contraction of $\beta_{2}$ 
over $\beta_{3}.$

\item Considering the parametric change of basis
\begin{equation*}
f_{t}(e_{1})=te_{1}, \quad
f_{t}(e_{2})=\frac{1}{2}(e_{1}+e_{2}),
\end{equation*} the limit when $t \longrightarrow 0$ gives an algebra 
isomorphic to $\beta_{4}.$

\item When considered in $\beta_1,$ the family of linear transformations
\begin{equation*}
f_{t}(e_{1})=\sqrt{\frac{t}{2}}(e_{1}+e_{2}), \quad
f_{t}(e_{2})=te_{2}
\end{equation*} defines an algebra isomorphic to $\beta_{5}$ in the limit.

\smallskip
The parametric change of basis
\begin{equation*}
f_{t}(e_{1})=te_{2}, \quad
f_{t}(e_{2})=t^{2}e_{1}
\end{equation*} induce a contraction of $\beta_{2}$ over $\beta_{5}.$

\smallskip
In the same way, if we take the family of non-singular endomorphisms
\begin{equation*} f_{t}(e_{1})=t(e_{1}+e_{2}) \quad
f_{t}(e_{2})=t^{2}e_{2},
\end{equation*} a contraction of $\beta_{3}$ over $\beta_{5}$ is obtained.

\smallskip
Finally, if we apply to $\beta_{4}$ the linear transformations 
\begin{equation*} f_{t}(e_{1})=e_{1}+te_{2},
\quad f_{t}(e_{2})=t^{2}e_{2},
\end{equation*} in the limit when $t\longrightarrow 0,$ it comes that $\beta_{5}$ is a contraction of $\beta_{4}.$\qedhere \end{enumerate}\end{proof}

\bigskip We summarize the classification and the contractions obtained for 
two dimensional real associative algebras in the following diagram, where 
each contraction is represented by an arrow. One can verify that 
rigid algebras correspond to those matrices which do not receive any arrow.

\bigskip

  \[ \xymatrix@C=1.5cm{*{\begin{pmatrix} 1 & 0 \\ 0 & 1 \\ 0 & 1 \\ -1 &
        0\end{pmatrix}} \ar[rd] \ar@/^-7ex/[rrdd] \ar@/^-10ex/[rrddd]
    & & *{\begin{pmatrix} 1 & 0 \\ 0 & 1 \\ 0 & 1 \\ 1 & 0
        \end{pmatrix}} \ar[ld] \ar[rd] \ar[dd] \ar@/_-7ex/[ddd]
      & & \\ & *{\begin{pmatrix} 1 & 0 \\ 0 & 1 \\ 0 & 1\\ 0 &
      0\end{pmatrix}} \ar[rd] \ar@/^-4ex/[rdd] & & *{\begin{pmatrix} 0
      & 0\\ 0 & 0\\ 0 & 0\\ 0 & 1 \end{pmatrix}} \ar[dl]
      \ar@/_-4ex/[ldd] & \\ *{\begin{pmatrix} 1 & 0 \\ 0 & 0\\ 0 & 1\\
      0 & 0\end{pmatrix}} \ar@/^-5ex/[rrd] & & *{\begin{pmatrix} 0 & 1
      \\ 0 & 0 \\ 0 & 0 \\ 0 & 0\end{pmatrix}} \ar[d] & &
      *{\begin{pmatrix} 1 & 0 \\ 0 & 1 \\ 0 & 0 \\ 0 & 0\end{pmatrix}}
      \ar@/_-5ex/[ldl] \\ & & *{\begin{pmatrix} 0 & 0\\ 0 & 0 \\ 0 & 0
      \\ 0 & 0
        \end{pmatrix}} & & } \]

{\bf Received: May 14, 2007}

\end{document}